\documentclass[12pt,reqno]{amsart}
\usepackage[utf8x]{inputenc}

\usepackage{latexsym,amssymb,amsmath,amsthm}
\usepackage{epsfig}
\usepackage{pst-eucl}
\usepackage{pst-plot,pstricks-add}

\newcommand{\R}{{\ensuremath{\mathbb{R}}}}
\newcommand{\N}{{\ensuremath{\mathbb{N}}}}

\renewcommand{\P}{\ensuremath{\mathbb{P}}}
\renewcommand{\dj}{d\kern-0.4em\char"16\kern-0.1em}
\newcommand{\E}{\ensuremath{\mathbb{E}}}

\newcommand{\cp}{\ensuremath{\mathrm{Cap}}}

\newtheorem{Thm}{Theorem}[section]

\newtheorem{Lem}[Thm]{Lemma}
\newtheorem{Prop}[Thm]{Proposition}

\theoremstyle{remark}

\theoremstyle{example}

\theoremstyle{definition}

\numberwithin{equation}{section}

\begin{document}
\bibliographystyle{amsalpha}

\title[Harnack inequality and H\" older regularity]{\bf Harnack inequality and H\" older regularity estimates for a L\' evy process with small jumps of high intensity}
\author{Ante Mimica}
\address{Department of Mathematics\\Bielefeld University\\Germany}
\email{amimica@math.hr}
\thanks{The author would like to thank Z. Vondra\v{c}ek and M. Kassmann for many helpful suggestions and comments.
}
\thanks{On leave from Department of Mathematics, University of Zagreb, Croatia.}
\thanks{Research supported in part by the MZOS Grant 037-0372790-2801 of the Republic of Croatia.}
\thanks{ Research supported in part by German Science Foundation DFG via IGK "Stochastics and real world models" and SFB 701.}
\date{}
\maketitle
\begin{abstract}
We consider a L\' evy process in $\R^d$ $ (d\geq 3)$ with the characteristic exponent
\[
\Phi(\xi)=\frac{|\xi|^2}{\ln(1+|\xi|^2)}-1.
\]
The scale invariant Harnack inequality and apriori estimates of harmonic functions in H\" older spaces  are proved.
\end{abstract}

\section{Introduction}

Let $X=(X_t,\P_x)_{t\geq 0,x\in\R^d}$ be a  L\' evy process in $\R^d$ ($d\geq 3$) with the characteristic exponent
\[
	\Phi(\xi)=\frac{|\xi|^2}{\ln(1+|\xi|^2)}-1.
\]

Let us give some motivation for the process $X$. It is known that the  variance gamma process can be obtained as a subordinate Brownian motion, where the corresponding subordinator is the gamma subordinator, i. e. a L\' evy process  whose Laplace exponent (cf. (\ref{eq:sub})) is given by
\[
	\ln(1+\lambda).
\] 
 It belongs to the class of the geometric stable processes (cf. \cite{SSV2}). The process $X$ is also a subordinate Brownian motion with subordinator that is a conjugate of the gamma subordinator. Namely, we take a subordinator with the Laplace exponent
\begin{equation}\label{eq:chexp1}
	\frac{\lambda}{\ln(1+\lambda)}-1.
\end{equation}
To avoid killing we subtract $1$ in (\ref{eq:chexp1}). Therefore we can say that  the process $X$ is (almost) conjugate  to the variance gamma process and so they are on the 'opposite sides'.

Another interesting property of  this process is that it is closer to the Brownian motion than any stable process. This can be argumented as follows.
Consider the potential operator $G$ defined by 
\[
	Gf(x)=\E_x\left[\int_0^\infty f(X_t)\,dt\right].
\]
This operator has a density, i.e. there exists a function $G(x,y)$, usually called the \emph{Green function}, with the following asymptotical properties (cf. Proposition \ref{prop:green_asymp}) 
\begin{equation}\label{eq:introt01}
	G(x,y)\sim \frac{\Gamma(\frac{d}{2}-1)}{4\pi^{d/2}}|x-y|^{2-d}\ln\frac{1}{|x-y|^2}\ \ \textrm{as }\ |x-y|\to 0+
\end{equation}
and
\begin{equation}\label{eq:introt02}
	G(x,y)\sim \frac{\Gamma(\frac{d}{2}-1)}{2\pi^{d/2}}|x-y|^{2-d}\ \ \textrm{as }\ |x-y|\to \infty.
\end{equation}
Comparing  (\ref{eq:introt01}) and (\ref{eq:introt02}) with  the Green function of the rotationally invariant $\alpha$-stable process  ($0<\alpha<2$): 
\[
	G^{(\alpha)}(x,y)=\frac{\Gamma(\frac{d}{2}-\frac{\alpha}{2})}{2^\alpha\pi^{d/2}\Gamma(\frac{\alpha}{2})}|x-y|^{\alpha-d}
\]
and the Green function of the Brownian motion in $\R^d$:
\[
	G^{(2)}(x,y)=\frac{\Gamma(\frac{d}{2}-1)}{4\pi^{d/2}}|x-y|^{2-d}.
\]
we see that the Green function of $X$ is "between" the Green functions of  Brownian motion and any stable process. We remark that $X$ is still a pure jump L\' evy process. 

The aim of this paper is to investigate some potential-theoretic notions of the process $X$. To be more precise, we investigate asymptotic behavior of the Green function and the L\' evy density of this process. Furthermore, we prove the scale invariant Harnack inequality and apriori H\" older estimates for the corresponding harmonic functions.

Although the origin of potential theory  is in the theory of differential equations,  it has also a probabilistic counterpart. The reason is that many local and non-local operators can be considered as infinitesimal generators of some Markov processes. 

Recently, probabilistic methods turned out to be very successful (cf. \cite{BL1}) in providing some steps in the proofs of the Harnack inequality and regularity estimates. Extensions to certain classes of L\' evy processes were obtained in \cite{SV1,BS,RSV,SSV2,KS,Mi}. More general jump processes were treated in \cite{BKa,BKa2,CK2}.

Before stating our results precisely, let us give a few comments on the main ingredient in the proof: a Krylov-Safonov type estimate. This kind of estimate for  jump processes appeared first in \cite{BL1}. In \cite{SV1} it was extended to some Markov processes. 

In these papers the following Krylov-Safonov type estimate was used:
\begin{equation}\label{eq:intro t05}
	\P_y(T_A<\tau_{B(0,r)})\geq c\,\frac{|A|}{|B(0,r)|},
\end{equation}
for all $r\in (0,1/2)$, $\ y\in B(0,r/2)$ and closed $A\subset B(0,r/2)$.
Here $T_A$ and $\tau_{B(x_0,r)}$ denote the first hitting time of $A$ and the first exit time from the ball $B(0,r)$, $|\cdot|$ denotes the Lebesgue measure in $\R^d$ and $c>0$ is a constant that does not depend on $r\in (0,1/2)$. 

Applying the same techniques to our case would lead to the  following estimate:
\begin{equation}\label{eq:intro t06}
	\P_y(T_A<\tau_{B(0,r)})\geq \frac{c}{\ln\frac{1}{r}}\,\frac{|A|}{|B(0,r)|}.
\end{equation}

The main difference between estimates  (\ref{eq:intro t05}) and (\ref{eq:intro t06}) is that the second one is not scale invariant. If we replace Lebesgue measure by some other set function, it is still possible  to obtain a scale invariant estimate of the similar type. More precisely, we have the following estimate (cf. Proposition \ref{heat:prop.1001})
\begin{equation}
\P_y(T_A<\tau_{B(x_0,r)})\geq c\,\frac{\cp(A)}{\cp\left(\overline{B(x_0,r)}\right)},
\end{equation}
where $\cp$ denotes the $0$-order capacity with respect to the process $X$ (cf. Section \ref{sect:prep}).

This idea appeared first in \cite{SSV2} and \cite{RSV}. We mention that in  the case of the geometric stable process considered in \cite{SSV2}, the known techniques lead to  the Krylov-Safonov type estimates that are not scale invariant neither with capacity nor with the Lebesgue measure.

In Section  \ref{sect:prep} we will see that the process $X$ is a purely discontinuous L\' evy process. Thus, the characteristic exponent $\Phi$ of $X$  is of the form
\begin{equation}\label{eq:chexp}
\Phi(\xi)=\int_{\R^d\setminus \{0\}}\left(e^{i\xi\cdot y}-1-i\xi\cdot y1_{\{|y|<1\}}\right)\Pi(dy).
\end{equation}
Here $\Pi$ denotes the \emph{L\' evy measure}, i.e. a measure on $\R^d\setminus \{0\}$ which satisfies the following integrability condition
\begin{equation}\label{eq:levycond}
\int_{\R^d\setminus\{0\}}(1\wedge |y|^2)\Pi(dy).
\end{equation}
Moreover, in our case $\Pi$ is absolutely continuous with respect to the Lebesgue measure:
\[
	\Pi(dy)=J(y)\,dy.
\]
The function $J$ is called  the \emph{L\' evy density}. At this point it is interesting to mention asymtpotical properties of $J$ (cf. Proposition \ref{prop:jump_asymp}):
\begin{equation}\label{eq:intt10}
	J(y)\sim \frac{4\Gamma(\frac{d}{2}+1)}{\pi^{d/2}}\cdot\frac{1}{|y|^{d+2}\left(\ln\frac{1}{|y|^2}\right)^2} \ \ \textrm{as }\ |y|\to 0.
\end{equation}
Comparing this with the L\' evy density of the rotationally invariant $\alpha$-stable process
\[
	J^{(\alpha)}(y)=\frac{\alpha2^{\alpha-1}\Gamma(\frac{d}{2}+\frac{\alpha}{2})}{\pi^{d/2}\Gamma(1-\frac{\alpha}{2})}\cdot\frac{1}{|y|^{d+\alpha}}
\]
we see that small jumps of the process $X$ are more intensive than the corresponding small jumps of any stable process. Using (\ref{eq:intt10}) we can see that the integrability condition of $\Pi$ given in (\ref{eq:levycond}) is barely satisfied.

We say that a function $h\colon \R^d\rightarrow [0,\infty)$ is \emph{harmonic} in an open set $D\subset \R^d$ with respect to the process $X$ if for any open set $B\subset D$ such that $B\subset \overline{B}\subset D$ the following is true
\[
	h(x)=\E_x[h(X_{\tau_B})1_{\{\tau_B<\infty\}}] \ \ \ \textrm{ for all }\  x\in B.
\]
Here $\tau_B=\inf\{t>0:X_t\not \in B\}$ is the \emph{first exit time} from $B$.
Denote by $B(x_0,r)$ the open ball in $\R^d$ with center $x_0\in \R^d$ and radius $r>0$.

The first result of this paper is the scale invariant Harnack inequality.
\begin{Thm}[Harnack inequality]\label{tm:harnack}
There exist $R>0$ and $L_1>0$ such that for any $x_0\in\R^d$ and $r\in (0,R)$ and any non-negative bounded function $h$ on $\R^d$ which is harmonic with respect to $X$ in $B(x_0,6r)$,
\[
 h(x)\leq L_1\,h(y)\ \textrm{ for all } x,y\in B(x_0,r).
\]
\end{Thm}

This type of Harnack inequality does not imply H\" older continuity directly via Moser's method of oscillation reduction. The relation of this two properties is currently investigated (cf. \cite{Kass}). 
Next result shows that harmonic functions locally satisfy uniform H\" older estimates.
\begin{Thm}[H\" older continuity]\label{tm:hoelder}
There exists $R'>0$, $\beta>0$ and $L_2>0$ such that for all $a\in \R^d$, $r\in (0,R')$ and any bounded function $h$ on $\R^d$ which is harmonic in $B(x_0,r)$ we have
\[
	|h(x)-h(y)|\leq L_2\|h\|_\infty r^{-\beta}|x-y|^\beta\ \textrm{ for all }\ x,y\in B(x_0,r/4).
\]
\end{Thm}

The paper is organized as follows. In Section \ref{sect:prep} we show that the process $X$ can be obtained as a subordinate Brownian motion and show some asymptotic properties of the L\' evy density and the Green function. We also prove a Krylov-Safonov type estimate. In Sections \ref{sect:hi} and \ref{sect:reg} we prove main results of the paper: Theorem \ref{tm:harnack} and Theorem \ref{tm:hoelder}.

\section{Preparatory results}\label{sect:prep}

A function $\phi\colon (0,\infty)\rightarrow (0,\infty)$ is called a \emph{Bernstein function}  if $\phi\in C^\infty((0,\infty))$ and 
\[
(-1)^{n-1} \phi^{(n)}(\lambda)\geq 0\ \textrm{ for all } \ n\in \N \ \textrm{ and }\  \lambda>0.
\]
We say that  $\phi\colon (0,\infty)\rightarrow (0,\infty)$ is a \emph{completely monotone function} if $\phi\in C^\infty((0,\infty))$ and
\[
(-1)^{n} \phi^{(n)}(\lambda)\geq 0\ \textrm{ for all } \ n\in \N\cup \{0\} \ \textrm{ and }\  \lambda>0.
\]

A \emph{subordinator} $S=(S_t)_{t\geq 0}$ is a L\' evy process taking values in $[0,\infty)$ and starting at 0. The Laplace transform of $S_t$ is given by
\[
	\E e^{-\lambda S_t}=e^{-t\phi(\lambda)},\ \lambda>0,
\]
where  $\phi$ is called the \emph{Laplace exponent} and  it is of the form (cf. \cite[p. 72]{Be})
\begin{equation}\label{eq:sub}
\phi(\lambda)=d\lambda+\int_{(0,\infty)}(1-e^{-\lambda t})\mu(dt).
\end{equation}
Here $d\geq 0$ and $\mu$ is a measure on $(0,\infty)$ (called a \emph{L\' evy measure})  satisfying
\[
\int_{(0,\infty)}(1\wedge t)\mu(dt)<\infty.
\]

Using \cite[Theorem 3.2]{SSV} we conclude that the Laplace exponent $\phi$ of $S$ is a Bernstein function. Conversely, if $\phi$ is a Bernstein function such that $\lim_{\lambda\to 0+}\phi(\lambda)=0$, then there exists a subordinator $S$ with the Laplace exponent $\phi$ (cf. \cite[Theorem I.1]{Be}).

We say that $f\colon (0,\infty)\rightarrow (0,\infty)$ is a \emph{ complete Bernstein function} if it has representation (\ref{eq:sub}) such that the L\' evy measure has a completely monotone density (with respect to the Lebesgue measure). It follows from \cite[Proposition 7.1]{SSV} that $f^*\colon (0,\infty)\rightarrow (0,\infty)$ defined by
\[
f^*(\lambda)=\frac{\lambda}{f(\lambda)}
\] is also a complete Bernstein function.

The \emph{potential measure} $U$ of the subordinator $S$ is defined by
\[
U(A)=\E \left[\int_0^\infty 1_{\{S_t\in A\}}\, dt\right],\ A\subset [0,\infty).
\]
The Laplace transform of $U$ is then
\begin{equation}\label{eq:lapl}
\mathcal{L}U(\lambda)=\int_{(0,\infty)}e^{-\lambda t}\,U(dt)=\frac{1}{\phi(\lambda)},\ \lambda>0.
\end{equation}

Let $B=(B_t,\P_x)_{t\geq 0, x\in \R^d}$ be a Brownian motion in $\R^d$ independent of the subordinator $S$. The process $X=(X_t,\P_x)_{t\geq 0,\ x\in \R^d}$ defined by
\[
X_t=B_{S_t},\ t\geq 0
\]
is called the \emph{subordinate Brownian motion}. By \cite[Theorem 30.1]{S} we conclude that $X$ is a L\' evy process and
\[
\E_x\left[e^{i\xi\cdot (X_t-X_0)}\right]=e^{-t\Phi(\xi)},\ \xi\in \R^d
\]
with $\Phi(\xi)=\phi(|\xi|^2)$. Moreover, we can rewrite $\Phi$ in the following way
\[
\Phi(\xi)=\int_{\R^d\setminus \{0\}}\left(e^{i\xi\cdot y}-1-i\xi\cdot y1_{\{|y|\leq 1\}}\right)j(|y|)\,dy
\]
where $j\colon (0,\infty)\rightarrow (0,\infty)$ is given by
\[
j(r)=(4\pi)^{-d/2}\int_{(0,\infty)}t^{-d/2}e^{-r^2/4t}\,\mu(dt),\ r>0.
\]
Therefore, the L\' evy measure of $X$ has density given by $J(y)=j(|y|)$.
Note that $j$ is a non-increasing function. 

From now on we denote by $S=(S_t)_{t\geq 0}$ the subordinator with the Laplace exponent
\[
\phi(\lambda)=\frac{\lambda}{\ln(1+\lambda)}-1
\]
and by $X$ the correspoding subordinate Brownian motion. 

Note that
\[
\ell(\lambda)=\ln(1+\lambda)=\int_0^\infty (1-e^{-\lambda t})\,t^{-1}e^{-t}\,dt
\]
and thus $\ell$ is a complete Bernstein function. Therefore, $\phi$ is also a complete Bernstein function and so the L\' evy measure of the subordinator $S$ has a completely monotone density $\mu(t)$. 

Let $T=(T_t)_{t\geq 0}$ be the subordinator with the Laplace exponent $\ell$. This process is known as the \emph{gamma subordinator}. 
 It follows from \cite[Corollary 10.7 and Corrolary 10.8]{SSV} that the potential measure $V$ of $T$ has a non-increasing density $v(t)$ and the following is true
 \begin{equation}\label{eq:connection}
 v(t)=1+\int_t^\infty \mu(s)\,ds,\ t>0. 
 \end{equation}
By \cite[Theorem 2.2]{SSV2} we get the following asymptotic behavior of $v$
\begin{equation}\label{eq:v_asym}
v(t)\sim t^{-1}\left(\ln\frac 1 t\right)^{-2}\ \textrm{ as } t\to 0+.
\end{equation}

Now we can prove the asymptotic behavior of the jumping function $J$. The proof of the following proposition is basically the proof of \cite[Lemma 3.1]{SSV2} but with the use of Potter's theorem (cf. \cite[Theorem 1.5.6 (ii)]{BGT}), which was also done in \cite[Lemma 5.1]{KSV2}.

\begin{Prop}\label{prop:jump_asymp}
The following asymptotic behavior of the function $j$ holds
  \[
   	j(r)\sim \frac{4\Gamma(\frac{d}{2}+1)}{\pi^{d/2}}\cdot \frac{1}{r^{d+2}\left(\ln\,\frac{1}{r^2}\right)^2}\ \textrm{ as } \ r\to 0+.
  \]
\end{Prop}
\proof
Using (\ref{eq:connection}) and (\ref{eq:v_asym}) we get
\[
 \int_t^\infty \mu(s)\,ds\sim \frac{1}{t(\ln{t})^2}\ \textrm{ as } \ t\to \,0+
\]
and thus by the Karamata's monotone density theorem (see \cite[Theorem 1.7.2]{BGT}) we have
\begin{equation}\label{eq:mu_asymp}
 \mu(t)\sim \frac{1}{t^2(\ln\,t)^2}\ \textrm{ as } \ t\to 0+.
\end{equation}
By change of variable we get
\begin{align}
j(r)&=(4\pi)^{-d/2}\int_0^\infty t^{-d/2}e^{-\frac{r^2}{4t}}\mu(t)\,dt\nonumber\\&=4^{-1}\pi^{-d/2}r^{-d+2}\int_0^\infty t^{d/2-2}e^{-t}\mu\left(\frac{r^2}{4t}\right)\,dt\nonumber\\
&=4^{-1}\pi^{-d/2}r^{-d+2}\mu(r^2)\int_0^\infty t^{d/2-2}e^{-t}\frac{\mu\left(\frac{|h|^2}{4t}\right)}{\mu(r^2)}\,dt.\label{eq:tmp101}
\end{align}
By Potter's theorem (cf.  \cite[Theorem 1.5.6 (iii)]{BGT}) we see that there is a constant $c_1>0$ such that 
\[
\frac{\mu\left(\frac{r^2}{4t}\right)}{\mu(r^2)}\leq c_1 (t^{2-1/2}\vee t^{2+1/2})\ \textrm{ for all } t>0 \textrm{ and } r>0.
\]
Therefore we can apply the dominated convergence theorem in (\ref{eq:tmp101}) and to get
\begin{equation}\label{eq:tmpj10}
\lim_{r\to 0+}\frac{j(r)}{4\pi^{-d/2}r^{-d+2}\mu(r^2)}=\Gamma\left(\frac{d}{2}+1\right).
\end{equation}
Combining (\ref{eq:mu_asymp}) and (\ref{eq:tmpj10}) we finish the proof. 
\qed

If $d\geq 3$, then by  \cite[Corollary 37.6]{S} and Proposition \ref{prop:jump_asymp} we conclude that $X$ is transient. Thus we can define a measure $G(x,\cdot)$ in the following way
\begin{equation}\label{eq:greenm}
G(x,A)=\E_x\left[\int_0^\infty 1_{\{X_t\in A\}}\,dt\right]=\int_0^\infty \P_x(X_t\in A)\,dt.
\end{equation}
Using \cite[Proposition 28.1]{S} we deduce that $X$ has a transition density $p(t,x,y)$ and the measure $G(x,\cdot)$ defined by (\ref{eq:greenm}) has a density which we denote by
\[
G(x,y)=\int_0^\infty p(t,x,y)\,dt
\]
and call the \emph{Green function} of $X$. Using \cite[Theorem 30.1]{S} we see that $G(x,y)=g(|x-y|)$, with
\begin{equation}\label{eq:greenj}
g(r)=(4\pi)^{-d/2}\int_{(0,\infty)} t^{-d/2}e^{-r^2/4t}\,U(dt),
\end{equation}
where $U$ is the potential measure of the subordinator $S$.
Combinig  \cite[Corollary 10.7 and Corollary 10.8]{SSV} we see that $U$ has a completely monotone density $u(t)$.

\begin{Lem}\label{prop:pot_asymp}
  We have the following asymptotics of $u$
\begin{align*}
 u(t)&\sim \ln\frac{1}{t}\ \textrm{ as } \ t\to 0+,&
u(t)&\to 2\ \textrm{ as } \ t\to \infty.
\end{align*}
\end{Lem}
\proof
We can readily check that
\begin{align*}
 	\phi(\lambda)&\sim \frac{\lambda}{2}\ \textrm{ as } \ \lambda\to 0+,&
 	\phi(\lambda)&\sim \frac{\lambda}{\ln\lambda}\ \textrm{ as } \ \lambda\to\infty,
\end{align*}
and thus by (\ref{eq:lapl}) we deduce 
\begin{align*}
 	\mathcal{L}U(\lambda)&\sim \frac{2}{\lambda}\ \textrm{ as } \ \lambda\to 0+, & 
 	\mathcal{L}U(\lambda)&\sim \frac{\ln\lambda}{\lambda}\ \textrm{ as } \ \lambda\to\infty.
\end{align*}
By the Karamata's Tauberian theorem (cf.  \cite[Theorem 1.7.1]{BGT}) we conclude that
\begin{align*}
 	U(t)&\sim 2t\ \textrm{ as } \ t\to \infty, &
 	U(t)&\sim t\ln\frac{1}{t}\ \textrm{ as } \ t\to 0+.
\end{align*}
Finally, using Karamata's monotone density theorem (cf. \cite[Theorem 1.7.2]{BGT}) we get
\begin{align*}
 	u(t)&\sim 2\ \textrm{ as } \ t\to \infty,
&
 	u(t)&\sim \ln\frac{1}{t}\ \textrm{ as } \ t\to 0+.
\end{align*}
\qed
\begin{Prop}\label{prop:green_asymp}
The following is true
	\begin{align*}
		g(r)&\sim \frac{\Gamma(\frac{d}{2}-1)}{2\pi^{d/2}}r^{2-d}\ln\frac{1}{r} \ \textrm{ as } \  r\to 0+, &
		g(r)&\sim \frac{\Gamma(\frac{d}{2}-1)}{2\pi^{d/2}}r^{2-d} \ \textrm{ as } \  r\to \infty.\\
	\end{align*}
\end{Prop}
\proof
Using (\ref{eq:greenj}) and changing variable we have
\begin{align}
g(r)&=4^{-1}\pi^{-d/2}r^{-d+2}\int_0^\infty t^{d/2-2}e^{-t}u\left(\frac{r^2}{4t}\right)\,dt\nonumber\\
&=4^{-1}\pi^{-d/2}r^{-d+2}u(r^2)\int_0^\infty t^{d/2-2}e^{-t}\frac{u\left(\frac{r^2}{4t}\right)}{u(r^2)}\,dt.\label{eq:tmp102}
\end{align}
From Potter's theorem (cf.  \cite[Theorem 1.5.6 (ii)]{BGT}) we deduce that there is a constant $c_1>0$ such that
\[
\frac{u\left(\frac{r^2}{4t}\right)}{u(r^2)}\leq c_1(t^{1/2}\vee t^{-1/2})\ \textrm{ for all } t>0 \textrm{ and } x\in \R^d,\ x\not= 0.
\]
Therefore we can use the dominated convergence theorem in (\ref{eq:tmp102}) to get
\begin{equation}\label{eq:tmpgj01}
\lim_{r\to 0+}\frac{g(r)}{4^{-1}\pi^{-d/2}r^{-d+2}u(r^2)}=\Gamma\left(\frac{d}{2}-1\right).
\end{equation}
Using (\ref{eq:tmpgj01}) and Lemma \ref{prop:pot_asymp} we obtain the asymptotics ti $0$. The other asymptotical formula is obtained simliarly.
\qed

We will need the following technical lemma later.
\begin{Lem}\label{heat:lem.100}
\begin{itemize}
\item[(a)]
 Let $f\colon (0,1)\rightarrow\R$ be defined by
\[
 f(t)=\frac{t^{d-2}}{\ln\frac 1 t}.
\]
Then $f$ is stricly increasing and 
\[
 f^{-1}(t)\sim (d-2)^{-\frac{1}{d-2}}\, t^{\frac{1}{d-2}}\left(\ln \frac 1 t\right)^{\frac{1}{d-2}}\ \textrm{ as } \ t\to 0+.
\]
\item[(b)] The following is true
\[
\int_0^rs^{d-1}g(s)\,ds \sim \frac{\Gamma(\frac{d}{2}-1)}{4\pi^{d/2}}r^2\ln\frac 1 r\ \textrm{ as }\ r\to 0+.
\]
\end{itemize}
\end{Lem}
\proof
\begin{itemize}
\item[(a)]
It is easy to see that $f$ is a strictly increasing function. Define $h\colon (0,1)\rightarrow \R$ by
\[
 h(t)=t^{\frac{1}{d-2}}\left(\ln \frac 1 t\right)^{\frac{1}{d-2}}.
\]
Then
\begin{align*}
 \lim_{t\to 0+}\frac{f^{-1}(t)}{h(t)}&=\lim_{t\to 0+}\frac{f^{-1}(f(t))}{h(f(t))}=\lim_{t\to 0+}\frac{t}{\frac{t}{\left(\ln\frac 1 t\right)^\frac{1}{d-2}}\left(\frac{t^{d-2}}{\ln\frac 1 t}\right)^\frac{1}{d-2}}\\
&=\lim_{t\to 0+}\left(\frac{\ln\frac 1 t}{(d-2)\ln\frac 1 t+\ln\ln\frac{1}{t}}\right)^\frac{1}{d-2}=(d-2)^{-\frac{1}{d-2}}.
\end{align*}
\item[(b)] By applying Karamata's theorem (cf. \cite[Proposition 1.5.8]{BGT}) we get
\begin{align*}
 \int_0^r s^{d-1}g(s)\,ds&\sim \frac{\Gamma(\frac{d}{2}-1)}{4\pi^{d/2}}\int_0^r s\ln\frac{1}{s^2}\,ds\\&\sim \frac{\Gamma(\frac{d}{2}-1)}{4\pi^{d/2}}r^2\ln \frac{1}{r}\ \textrm{ as } \ r\to 0+.
 \end{align*}
 \end{itemize}
\qed

Let $D\subset \R^d$ be an open set. We define the \emph{killed process} $X^D=(X^D_t)_{t\geq 0}$ by killing process $X$ upon exiting set $D$, i.e.
\[
X^D_t=\left\{
\begin{array}{cl}
X_t, & t<\tau_D\\
\partial,& t\geq \tau_D.
\end{array}
\right.
\] 
Here $\partial$ is an extra point adjoined to $D$. In this case the killed process also has a transition density and it is given (cf. proof of \cite[Theorem 2.4]{CZ}) by
\begin{equation}\label{eq:tranD}
p_D(t,x,y)=p(t,x,y)-\E_x \left[p(t-\tau_D,X_{\tau_D},y);\tau_D<t\right].
\end{equation}
The Green function of $X^D$ also exists and it is given by
\begin{equation}\label{eq:greenD}
G_D(x,y)=\int_0^\infty p_D(t,x,y)\,dt=G(x,y)-\E_x\left[G(X_{\tau_D},y)\right]\ \textrm{ for } x,y\in D.
\end{equation}
 Since $\P_x(X_{\tau_{B(x_0,r)}}\in \partial B(x_0,r))=0$ (cf. \cite{SZ}) for $x_0\in \R^d$ and $r>0$, it follows from Theorem 1 in \cite{IW} that for any non-negative function $h\colon\R^d\rightarrow [0,\infty)$ we have
 \begin{equation}\label{eq:tmp120}
 \E_x[h(X_{\tau_{B(x_0,r)}})]=\int_{\overline{B(x_0,r)}^c}\int_{B(x_0,r)} G_{B(x_0,r)}(x,y)j(|z-y|)h(z)\,dy\,dz.
 \end{equation}
 If we define \emph{Poisson kernel} $K_{B(x_0,r)}\colon B(x_0,r)\times \overline{B(x_0,r)}^c\rightarrow [0,\infty)$ by
  \[
 K_{B(x_0,r)}(x,z)=\int_{B(x_0,r)}G_{B(x_0,r)}(x,y)j(|z-y|)\,dy\ \textrm{ for } x\in B(x_0,r),\ z\in \overline{B(x_0,r)}^c,
 \]
from (\ref{eq:tmp120}) we get 
 \begin{equation}\label{eq:Poisson1}
 \E_x[h(X_{\tau_{B(x_0,r)}})]=\int_{\overline{B(x_0,r)}^c} K_{B(x_0,r)}(x,z)h(z)\,dz.
 \end{equation}
 
 \begin{Prop}\label{heat:prop.80}
 There exist constants  $R_0\in(0,1/6)$ and $C_1>0$ such that for  any  $r\leq R_0$ and $x_0\in\R^d$,
\begin{equation}
G_{B(x_0,4r)}(x,y)\geq C_1r^{2-d}\ln\frac{1}{r}\  \textrm{ for all }  x,y\in B(x_0,r).
\end{equation}
\end{Prop}
\proof
Choose $0<c_1<1<c_2$ such that \[c_1^2\left(\frac{1}{2}\right)^{d-2}-c_2^2\left(\frac{1}{3}\right)^{d-2}> 0.\]
Using Proposition \ref{prop:green_asymp}  we can choose $R_0\in(0,1/6)$ such that for $r\leq 3R_0$ we have
\begin{equation}\label{est_1}
 c_1c_3r^{2-d}\ln\frac{1}{r}\leq g(r)\leq c_2c_3r^{2-d}\ln\frac{1}{r},\ \ \ c_1\leq \frac{\ln\frac{1}{2r}}{\ln\frac{1}{r}}\leq c_2,\ \ \  c_1\leq \frac{\ln\frac{1}{3r}}{\ln\frac{1}{r}}\leq c_2,
\end{equation}
where $c_3=\frac{\Gamma(d/2-1)}{2\pi^{d/2}}$.
Let $r\leq R_0$, $x_0\in\R^d$ and  $x,y\in B(x_0,r)$. By (\ref{est_1}) and monotonicity of $g$ we get
\begin{align*}
 G_{B(x_0,4r)}(x,y)&=G(x,y)-\E_x[G(Y_{\tau_{B(x_0,4r)}},y)]=g(|x-y|)-\E_x[g(|Y_{\tau_{B(x_0,4r)}}-y|)]\\
&\geq g(2r)-g(3r)\geq c_3\left( c_1(2r)^{2-d}\ln\frac{1}{2r}-c_2\left(3r\right)^{2-d}\ln\frac{1}{3r}\right)\\
&=c_3r^{2-d}\ln{\frac{1}{r}}\left(c_1\,\left(\frac{1}{2}\right)^{d-2}\frac{\ln{\frac{1}{2r}}}{\ln{\frac{1}{r}}}-c_2\,\left(\frac{1}{3}\right)^{d-2}\frac{\ln{\frac{1}{3r}}}{\ln{\frac{1}{r}}}\right)\\&\geq c_3r^{2-d}\ln{\frac{1}{r}}\left(c_1^2\left(\frac{1}{2}\right)^{d-2}-c_2^2\left(\frac{1}{3}\right)^{d-2}\right).
\end{align*}
Hence we may take
\[C_1=c_3\left(c_1^2\left(\frac{1}{2}\right)^{d-2}-c_2^2\left(\frac{1}{3}\right)^{d-2}\right)>0.\]
\qed

\begin{Prop}\label{prop:Poisson1}
 There exist $R_1\in(0,R_0]$ and  a constant $C_2>0$ such that for  any $r\leq R_1$ and $x_0\in\R^d$,
\begin{equation}
 K_{B(x_0,r)}(x,z)\leq C_2\,K_{B(x_0,4r)}(y,z)\ \textrm{ for all } x,y\in B(x_0,r/2),\,z\in B(x_0,4r)^c.
\end{equation}
\end{Prop}
\proof
Take $r\leq R_0$, $x_0\in\R^d$, $x,y\in B(x_0,r/2)$ and $z\in B(a,3r)^c$. Using Proposition \ref{heat:prop.80} we get
\begin{align}\nonumber
 K_{B(x_0,4r)}(y,z)&=\int_{B(x_0,4r)}G_{B(x_0,4r)}(y,u)\,j(|z-u|)\,du\\&\geq\int_{B(x_0,r)}G_{B(x_0,4r)}(y,u)\,j(|z-u|)\,du\\\label{heat:eq.34}
&\geq C_1r^{2-d}\ln\frac{1}{r}\int_{B(x_0,r)}j(|z-u|)\,du.
\end{align}
On the other side, applying  \cite[Lemma 2.7]{Mi} to $j$, and then using Proposition \ref{prop:green_asymp}
 and Lemma \ref{heat:lem.100} (b)  we see that there exist constants  $R_1\in(0,R_0]$ and 
$c_1,c_2,c_3>0$ such that\\ 
\begin{align*}
 K_{B(x_0,r)}(x,z)=&\,\int_{B(x_0,3r/4)}G_{B(x_0,r)}(x,v)j(|z-v|)\,dv\\&+\int_{B(a,
r)\setminus B(x_0,3r/4)}G_{B(x_0,r)}(x,v)j(|z-v|)\,dv\\
\leq\,&
c_1r^{-d}\int_{B(x_0,r)}j(|z-u|)\,du\int_{B{(a,3r/4)}}g(|x-v|)\,dv\,\\&+c_2r^{2-d}
\ln\frac{1}{r}\int_{B(x_0,r)\setminus B(x_0,3r/4)}j(|z-u|)\,du\\
\leq\,& c_3\,r^{2-d}\ln\frac{1}{r}\int_{B(x_0,r)}j(|z-u|)\,du+c_2r^{2-d}\ln\frac{1}{r}\int_{B(x_0,r)}j(|z-u|)\,du\\
=&(c_2+c_3)r^{2-d}\ln\frac{1}{r}\int_{B(x_0,r)}j(|z-u|)\,du,
\end{align*}
where in the first term in the last inequality we have used
\[
 \int_{B(x_0,3r/4)}g(|x-v|)\,dv\leq \int_{B(x,2r)}g(|x-v|)\,dv\leq c_1'r^2\ln\frac{1}{r}.
\]
Therefore
\[
 K_{B(x_0,r)}(x,z)\leq \frac{c_2+c_3}{C_1}K_{B(x_0,4r)}(y,z).
\]
\qed

 For a measure $\rho$ on $\R^d$ we define its \emph{potential} by
 \[
 G\rho(x)=\int_{\R^d}G(x,y)\,\rho(dx).
 \]
Denote by  $\cp$ the (0-order) \emph{capacity} with respect to $X$ (cf.  \cite[Section II.2]{Be}). It is proved in \cite[Corollary II.8]{Be} that for any compact set $K\subset \R^d$ there exists a measure $\rho_K$, called the  \emph{equilibrium measure}, which is supported by $K$ and satisfies 
\begin{equation}\label{eq:equilibrium1}
G\rho_K(x)=\P_x(T_K<\infty)\ \textrm{ for a.e. } x\in \R^d.
\end{equation}
Moreover, the following is true
\[
\frac{1}{\cp (K)}=\inf\left\{\int_{\R^d} G\rho(x)\,\rho(dx)\colon \rho \textrm{ is a probability measure supported by}\ K\right\}
\]
and the infimum is attained at the equilibrium measure $\rho_K$. If we combine Lemma \ref{heat:lem.100} (b) and  \cite[Proposition 5.3]{SSV2} we conclude that there exist constants $C_3, C_{4}>0$ such that 
\begin{equation}\label{eq:est_cap}
C_3\frac{r^{d-2}}{\ln\frac{1}{r}} \leq \cp\left(\overline{B(x_0,r)}\right)\leq C_{4}\frac{r^{d-2}}{\ln\frac{1}{r}}\ \textrm{ for } x_0\in\R^d, \ 0<r\leq 1/2.
\end{equation}
Now we can prove a Krylov-Safonov-type estimate.
\begin{Prop}\label{heat:prop.1001}
 There exists a constant $C_{5}>0$  such that for any $x_0\in\R^d$, $r\leq R_1$, closed subset A of $B(x_0,r)$ and $y\in B(x_0,r)$,
\begin{equation*}
 \P_y(T_A<\tau_{B(x_0,4r)})\geq C_{5}\frac{\cp(A)}{\cp\left(\overline{B(x_0,4r)}\right)}.
\end{equation*}
\end{Prop}
\proof
Let $x_0\in\R^d$, $r\leq R_1$ and let $A\subset B(x_0,r)$ be a closed subset. We may assume that $\cp(A)>0$. Let $\rho_A$ be the equilibrium measure of $A$. If $G_{B(x_0,4r)}$ is the Green function of the process $X$ killed upon exiting from $B(x_0,4r)$, then  for $y\in B(x_0,r)$ we have
\begin{align}
 G_{B(x_0,4r)}\rho_A(y)&=\E_y[G_{B(x_0,4r)}\rho_A(Y_{T_A});T_A<\tau_{B(x_0,4r)}]\nonumber\\
&\leq \sup_{z\in\R^d}G_{B(x_0,4r)}\rho_A(z)\P_y(T_A<\tau_{B(x_0,4r)})\nonumber\\
&\leq \P_y(T_A<\tau_{B(x_0,4r)}),\label{heat:eq.40}
\end{align}
since $G_{B(x_0,4r)}\rho_A(z)\leq G\rho_A(z)\leq 1$ by (\ref{eq:equilibrium1}). 
Also, for any $y\in B(x_0,r)$ we have
\begin{align}
 G_{B(x_0,4r)}\rho_A(y)&=\int_{\R^d}G_{B(x_0,4r)}(y,z)\rho_A(dz)\geq
\rho_A(\R^d)\inf_{z\in B(x_0,r)}G_{B(x_0,4r)}(y,z)\nonumber\\
&=\cp(A)\inf_{z\in B(x_0,r)}G_{B(x_0,4r)}(y,z).\label{heat:eq.41}
\end{align}
Using (\ref{heat:eq.40}), (\ref{heat:eq.41}) and Proposition \ref{heat:prop.80} we  obtain
\[
 \P_y(T_A<\tau_{B(x_0,4r)})\geq C_1\cp(A)r^{2-d}\ln\frac{1}{r}.
\]
By (\ref{eq:est_cap}) we see that
\[
 \P_y(T_A<\tau_{B(x_0,4r)})\geq C_1C_3\frac{\cp(A)}{\cp\left(\overline{B(x_0,4r)}\right)}
\]
for $r\leq R_1$.
\qed

\section{Harnack inequality}\label{sect:hi}

{\bf Proof of Theorem \ref{tm:harnack}.} Define $f(t)= \frac{t^{d-2}}{\ln\frac{1}{t}}$. By Lemma \ref{heat:lem.100} (a) we can choose $R\leq \frac{R_1}{4}\wedge\frac{1}{16}$ such that 
\begin{equation}\label{heat:tmp}
 f^{-1}(r)\leq c_0r^{\frac{1}{d-2}}\left(\ln\frac 1 r\right)^{\frac{1}{d-2}}\ \textrm{ for all } \ r\leq R
\end{equation}
and for some constant $c_0\geq 1$.

Let $x_0\in\R^d$ and  $r\leq R$.
	Without loss of generality we may suppose
	\[
	 	\inf_{z\in B(x_0,r)} h(z)=\frac{1}{2}
	\]
Let $z_0\in B(x_0,r)$ be such that $h(z_0)\leq 1$. It is enough to show that $h$ is bounded from above by some constant independent of $h$. By Proposition \ref{heat:prop.1001}  there exists $c_1>0$ such that
\begin{equation}\label{eq:har1}
 	\P_x(T_F<\tau_{B(x,s)})\geq c_1,
\end{equation}
for any $s\in (0,R_1)$, $x\in \R^d$ and a compact subset $F\subset B(x,s/4)$ such that
\[
 \frac{\cp(F)}{\cp\left(\overline{B(x,s/4)}\right)}\geq \frac{1}{3}.
\]
Put
\begin{equation}\label{eq:eta}
 \eta=\frac{c_1}{3},\ \ \zeta=\frac{\eta}{3}\wedge \frac{\eta}{C_2}.
\end{equation}
Suppose that there exists $x\in B(x_0,r)$ such that $h(x)=K$ for
\begin{equation}\label{eq:K_1}
 K>c_0^{d-1}4^{3(d-2)}\,\left(\frac{2C_4}{C_3C_5\zeta}\right)^{1/2}.
\end{equation}
It is possible to choose a unique $s>0$ such that
\[
f(s/4)=\frac{2\,C_4}{C_3\,C_5\,\zeta\,K}f(4r),
\]
since $f$ is strictly increasing and continuous on $(0,1)$ and $c_2:=\frac{2\,C_4}{C_3\,C_5\,\zeta\,K}<\frac{1}{4^3}$.

Using inequality 
\[
 p\ln\frac 1 p<2\sqrt{p}\ \textrm{ for } \ p\in (0,1)
\]
and (\ref{heat:tmp}) we get
\begin{align}
 \frac s 4&\leq c_0\left(c_2f(4r)\right)^{\frac{1}{d-2}}\left(\ln\frac{1}{c_2f(4r)}\right)^{\frac{1}{d-2}}\nonumber\\
 &\leq c_0^{1+\frac{1}{d-2}}c_2^{\frac{1}{d-2}}\left(\frac{\ln\frac{1}{c_2}}{\ln\frac{1}{4r}}+d-2-\frac{\ln\ln\frac{1}{4r}}{\ln\frac{1}{4r}}\right)^{\frac{1}{d-2}}4r\nonumber\\
 &\leq  c_0^{1+\frac{1}{d-2}}\left(c_2\ln\frac{1}{c_2}+(d-2)c_2\right)^{\frac{1}{d-2}}4r\nonumber\\
 &\leq  c_0^{1+\frac{1}{d-2}}\left(2\sqrt{c_2}+(d-2)c_2\right)^{\frac{1}{d-2}}4r\nonumber\\
 &\leq c_0^{1+\frac{1}{d-2}}(d\sqrt{c_2})^{\frac{1}{d-2}}4r\label{eq:tmpharn01}
 \end{align}
and thus $s\leq \frac{r}{4}$ by (\ref{eq:K_1}).

By (\ref{eq:est_cap}) we obtain
\[
 \cp\left(\overline{B(x,s/4)}\right)\geq C_3f(s/4)\geq \frac{2}{C_5\,\zeta\,K}\cp\left(\overline{B(x_0,4r)}\right).
\]
Let $A$ be a compact subset of
\[
 A'=\{t\in B(x,s/4)\colon h(t)\geq \zeta K \}.
\]
By the optional stopping theorem  we have
\begin{align}
 1&\geq h(z_0)\geq E_{z_0}[h({X_{T_A\wedge \tau_{B(x_0,4r)}}});T_A< \tau_{B(x_0,4r)}]\nonumber\\
&\geq \zeta K\P_{z_0}(T_A< \tau_{B(x_0,4r)})\geq C_5\zeta K\frac{\cp(A)}{\cp\left(\overline{B(x_0,4r)}\right)},\label{eq:tmpharn02}
\end{align}
where in the last inequality we have used Proposition \ref{heat:prop.1001}. Therefore, from (\ref{eq:tmpharn01}) and (\ref{eq:tmpharn02}) we conclude 
\[
 \frac{\cp(A')}{\cp\left(\overline{B(x,s/4)}\right)}= \frac{\cp(A')}{\cp\left(\overline{B(x_0,4r)}\right)}\cdot  \frac{\cp\left(\overline{B(x_0,4r)}\right)}{\cp\left(\overline{B(x,s/4)}\right)}\leq \frac{1}{2}
\]
and thus, by subadditivity of capacity, there exists a compact set \[F\subset \overline{B(x,s/4)}\setminus A'\] such that
\begin{equation}\label{eq:f}
 \frac{\cp(F)}{\cp\left(\overline{B(x,s/4)}\right)}\geq \frac{1}{3}.
\end{equation}
 Next we prove that
\[
 \E_x[h(X_{\tau_{B(x,s)}}); X_{\tau_{B(x,s)}}\not\in B(x,4s)]\leq \eta K.
\]
If the latter is not true, then \[\E_x[h(X_{\tau_{B(x,s)}}); X_{\tau_{B(x,s)}}\not\in B(x,4s)]> \eta K\] and by Proposition \ref{prop:Poisson1} and (\ref{eq:Poisson1}) for any $y\in B(x,s/4)$ we have
\begin{align*}
 h(y)&=\E_y[h(X_{\tau_{B(x,4s)}})]=\E_y[h(X_{\tau_{B(x,4s)}});X_{\tau_{B(x,4s)}}\not\in B(x,4s)]\\
&=\int_{\overline{B(x,4s)}^c} K_{B(x,4s)}(y,z)h(z)\,dz\geq C_2^{-1}\int_{\overline{B(x,4s)}^c} K_{B(x,s)}(y,z)h(z)\,dz\\
&=C_2^{-1}\E_y[h(X_{\tau_{B(x,s)}});X_{\tau_{B(x,s)}}\not\in B(x,4s)]\geq C_2^{-1}\eta K\geq \zeta K,
\end{align*}
which is a contradiction with (\ref{eq:f}) and the definition of the set $A'$.

Set
\[
 M=\sup_{B(x,4s)}\,h.
\]
We have
\begin{align*}
 K&=h(x)=\E_x[h(Y_{\tau_{B(x,s)}})]\\
&=\E_x[h(Y_{T_F});T_F<\tau_{B(x,s)}]+\E^x[h(Y_{\tau_{B(x,s)}});\tau_{B(x,s)}<T_F,X_{\tau_{B(x,s)}}\in B(x,4s)]\\
&+\E_x[h(Y_{\tau_{B(x,s)}});\tau_{B(x,s)}<T_F,X_{\tau_{B(x,s)}}\not\in B(x,4s)]\\
&\leq \zeta K\P_x(\tau_{B(x,s)}<T_F)+M\P_x(\tau_{B(x,s)}<T_F)+\eta K\\
&= \zeta K\P_x(\tau_{B(x,s)}<T_F)+M(1-\P_x(\tau_{B(x,s)}<T_F))+\eta K
\end{align*}
and thus
\[
 \frac{M}{K}\geq \frac{1-\eta-\zeta \P^x(\tau_{B(x,s)}<T_F)}{1-\P^x(\tau_{B(x,s)}<T_F)}\geq 1+2\beta,
\]
for some $\beta>0$. It follows that there exists $x'\in B(x,4s)$ such that $h(x')\geq K(1+\beta)$. Repeating this procedure, we get a sequence $(x_n)$ such that $h(x_n)\geq K(1+\beta)^{n-1}$ and 
\[
|x_{n+1}-x_n|\leq \left(4^3c_0\left(dc_0\sqrt{c_2}\right)^{\frac{1}{d-2}}\right)^nr=\left(c_3K^{-\frac{1}{d-2}}\right)^n\,r
\]
Therefore,
\[
 \sum_{n=1}^\infty |x_{n+1}-x_n|\leq c_4K^{-\frac{1}{d-2}}r.
\]

If $K>c_4^{d-2}$ we can find sequence $(x_n)$ in $B(x_0,2r)$ such that $h(x_n)\to \infty$ which is a contradiction to $h$ being bounded and so
\[
 \sup_{x\in B(x_0,r)}h(x)\leq c_4^{d-2}.
\]
\qed

\section{Regularity}\label{sect:reg}
 \begin{Lem}\label{lem:tmp200}
 There exists a  constant $C_6>0$ such that for all $r\in (0,1/8)$, $s\in [4r,1/2)$ and $x_0\in \R^d$ we have
 \[
 \P_x(X_{\tau_{B(x_0,r)}}\not \in B(x_0,s))\leq C_6\frac{r^2\ln\frac 1 r}{s^2\left(\ln\frac 1 s\right)^2}\ \textrm{ for all }\ x\in B(x_0,r/2).
 \]
 \end{Lem}
 \proof Let $r\in (0,1/8)$, $s\in [4r,1/2)$ and $x_0\in \R^d$. 
 Using (\ref{eq:tmp120}) with $h=1_{B(x_0,s)^c}$ for $x\in B(x_0,r/2)$ we have
 \begin{align*}
 \P_x(X_{\tau_{B(x_0,r)}}\not \in B(x_0,s))&=\int_{B(x_0,r)}G_{B(x_0,r)}(x,u)\int_{B(x_0,s)^c}j(|z-u|)\,dz\,du\\
 &\leq \int_{B(x_0,r)}g(|u-x|)\int_{B(u,s/2)^c}j(|z-u|)\,dz\,du\\
 &=\int_{B(x-x_0,r)}g(|u|)\,du\cdot\int_{B(0,s/2)^c}j(|z|)\,dz\\
 &\leq \int_{B(0,2r)}g(|u|)\,du\cdot\int_{B(0,s/2)^c}j(|z|)\,dz,
 \end{align*}
 where in the second inequality we have used the facts that $B(u,s/2)\subset B(x_0,s)$ and $G_{B(x_0,r)}(y,u)\leq g(|u-y|)$, while in the last inequality we have used $B(x-x_0,r)\subset B(0,2r)$. Now the conclusion follows from Proposition \ref{prop:jump_asymp} and Proposition \ref{prop:green_asymp}. 
 \qed
 
 {\bf Proof of Theorem \ref{tm:hoelder}.} Let $r\in (0,1/6)$, $x_0\in \R^d$ and let  $h\colon \R^d\rightarrow [0,\infty)$ be bounded by $M>0$ and harmonic in $B(x_0,r)$. 
 
 Let $z_0\in B(x_0,r/4)$. Define
 \[
 	r_n=\gamma_1 \,4^{-n},
 \]
 where we choose $\gamma_1>0$ small enough so that $B(x_0,4r_1)\subset B(z_0,r/4)$. Set $B_n=B(z_0,r_n)$ and $\tau_n=\tau_{B_n}$ and 
 \[
 s_n= b^{-n},
 \]
 where constant $b>1$ will be chosen later. 
 Let
 \[
 m_n=\inf_{x\in B_n}h(x)\ \ \ \textrm{ and }\ \ \  M_n=\sup_{x\in B_n}h(x).
 \]
 It is enough to prove that 
 \begin{equation}\label{eq:ind}
 M_k-m_k\leq s_k
 \end{equation}
 for all $k\geq n_0$ for some $n_0\in \N$.
 We prove this by induction. Assume that (\ref{eq:ind}) holds for $k\in \{n_0,n_0+1,\ldots,n\}$, where $n\geq n_0$.
 
 Let $\varepsilon>0$ and let $x,y\in B_{n+1}$ such that $h(x)\leq m_{n+1}+\varepsilon$ and $h(y)\geq M_{n+1}-\varepsilon$. Our aim is to show that $h(y)-h(x)\leq s_{n+1}$. Then we have
 \[
 M_{n+1}-m_{n+1}\leq s_{n+1}+2\varepsilon
 \]
 and since $\varepsilon>0$ is arbitrary, we get (\ref{eq:ind}) for $k=n+1$.
 
 Let $A=\left\{x\in B_{n+1}\colon h(x)\leq \frac{m_n+M_n}{2}\right\}$ and assume that 
 \[
 	\frac{\cp(A)}{\cp(B_{n+1})}\geq \frac{1}{2}
 \]
 (if this is not true, then we consider function $M-h$ and use the subadditivity of capacity). By Choquet's theorem $A$ is capacitable and therefore there exists a compact subset $K\subset A$ such that 
 \[
	 \frac{\cp(K)}{\cp(B_{n+1})}\geq \frac{1}{3}.
 \]

 By the optional stopping theorem, we have
 \begin{align*}
 h(x)-h(y)=&\ \E_x[h(X_{\tau_n\wedge T_K})-h(y)]\\
 =&\ \E_x[h(X_{\tau_n\wedge T_K})-h(y);T_K<\tau_n,X_{\tau_n}\in B_{n-1}\setminus B_n]\\ &+\E_x[h(X_{\tau_n\wedge T_K})-h(y); T_K>\tau_n, X_{\tau_n}\in B_{n-1}\setminus B_n]\\
 &+\sum_{i=1}^{n-2}\E_x[h(X_{\tau_n\wedge T_K})-h(y);X_{\tau_n}\in B_{n-i-1}\setminus B_{n-i}]\\
 &+\E_x[h(X_{\tau_n\wedge T_K})-h(y); X_{\tau_n}\not\in B_1]\\
 \leq & \left(\frac{m_n+M_n}{2}-m_n\right)\P_x(T_K<\tau_n)+(M_{n-1}-m_{n-1})\P_x(T_K>\tau_n)\\
 &\ + \sum_{i=1}^{n-2}(M_{n-i-1}-m_{n-i-1})\P_x(X_{\tau_n}\not \in B_{n-i-1})+2M\P_x(X_{\tau_n}\not\in B_1).
  \end{align*}
  
  It follows from Proposition \ref{heat:prop.1001} that there is a constant $c_1>0$ such that 
  \begin{equation}\label{eq:tmphold01}
  p_n:=\P_x(T_K<\tau_n)\geq c_1.
  \end{equation}
By Lemma \ref{lem:tmp200} we obtain
\[
h(x)-h(y)\leq \frac{1}{2}s_n\,p_n+s_{n-1}(1-p_n)+C_6\sum_{i=1}^{n-2} s_{n-i-1}\frac{r_n^2\ln\frac{1}{r_n}}{r_{n-i}^2\left(\ln\frac{1}{r_{n-i}}\right)^2}+2MC_6\frac{r_n^2\ln\frac{1}{r_n}}{r_{1}^2\left(\ln\frac{1}{r_{1}}\right)^2}.
\]
Then there exist constants $c_2,c_3>0$ such that
\begin{align}\nonumber 
&h(x)-h(y)\leq\nonumber\\&\leq s_{n+1}\left[\frac{b}{2}\,p_n+b^2(1-p_n)+c_2\,n\left(4^{2}b^{-1}\right)^{-n}\,b^2\sum_{i=2}^{n-1}\frac{(4^2b^{-1})^{-i}}{i^2}+c_3Mb \,n\,(4^2b^{-1})^{-n}\right]\nonumber\\
&\leq s_{n+1}\left[\frac{b}{2}\,p_n+b^2(1-p_n)+c_2\,n\left(4^{2}b^{-1}\right)^{-n}\,b^2\sum_{i=2}^{\infty}\frac{(4^2b^{-1})^{-i}}{i^2}+c_3Mb \,n\,(4^2b^{-1})^{-n}\right].\label{eq:tmp234}
\end{align}
Choose  $b\in (1,2)$ such that $(1-c_1)b^2+c_1\,\frac{b}{2}<1$. By (\ref{eq:tmphold01}) we conclude 
\[
	\frac{1}{2}b\,p_n+b^2(1-p_n)=b^2-p_n\left(b^2-\frac{b}{2}\right)\leq b^2-c_1\left(b^2-\frac{b}{2}\right)<1.
\]
Since $4^2b^{-1}>1$, we see that the last two terms in the parenthesis in (\ref{eq:tmp234}) can be made arbitrary small for $n$ large enough. Thus there is $n_0\in\N$ such that for any $n\geq n_0$ we have
\[
h(x)-h(y)\leq s_{n+1},
\]
which was to be proved.
\qed

\textbf{Acknowledgements.} The author would like to thank Z. Vondra\v{c}ek and M. Kassmann for many helpful suggestions and comments.

\end{document}